\documentclass [reqno]{amsart}
\usepackage {latexsym,epsfig}
\usepackage [latin1]{inputenc}

\newcommand {\PP}{{\mathbb P}}
\newcommand {\QQ}{{\mathbb Q}}
\newcommand {\ZZ}{{\mathbb Z}}
\newcommand {\calC}{{\mathcal C}}
\newcommand {\calO}{{\mathcal O}}
\newcommand {\calP}{{\mathcal P}}
\newcommand {\calX}{{\mathcal X}}

\newcommand {\preprint}[2]{preprint \discretionary {#1/}{#2}{#1/#2}}

\newtheorem {theorem}{Theorem}[section]
\newtheorem {lemma}[theorem]{Lemma}
\newtheorem {proposition}[theorem]{Proposition}
\newtheorem {corollary}[theorem]{Corollary}
\newtheorem {definition}[theorem]{Definition}
\theoremstyle {definition}
\newtheorem {example}[theorem]{Example}
\theoremstyle {remark}
\newtheorem {remark}[theorem]{Remark}

\begin {document}


\title [The number of plane conics 5-fold tangent to a given curve]{%
  The number of plane conics 5-fold tangent to a given curve}
\author {Andreas Gathmann}
\address {Institute for Advanced Study, School of Mathematics,
          1 Einstein Drive, Princeton NJ 08540, USA}
\email {andreas@ias.edu}

\begin {abstract}
  Given a general plane curve $Y$ of degree $d$, we compute the number $ n_d $
  of irreducible plane conics that are 5-fold tangent to $Y$. This problem has
  been studied before by Vainsencher \cite {V} using classical methods, but it
  could not be solved there because the calculations received too many
  non-enumerative correction terms that could not be analyzed. In our current
  approach, we express the number $ n_d $ in terms of relative Gromov-Witten
  invariants that can then be directly computed. As an application, we consider
  the K3 surface given as the double cover of $ \PP^2 $ branched along a sextic
  curve. We compute the number of rational curves in this K3 surface in the
  homology class that is the pull-back of conics in $ \PP^2 $, and compare this
  number to the corresponding Yau-Zaslow K3 invariant. This gives an example of
  such a K3 invariant for a non-primitive homology class.
\end {abstract}

\maketitle


Let $ Y \subset \PP^2 $ be a generic plane curve of degree $ d \ge 5 $. We want
to consider smooth plane conics that are 5-fold tangent to $Y$. As the space of
all plane conics is 5-dimensional and each tangency imposes one condition on
the curves, we expect a finite number of such 5-fold tangent conics. It will be
easy to see that this number is indeed finite; let us call it $ n_d $. The goal
of this paper is to compute it.

Of course this is a classical problem, and attempts have been made to solve it
using classical methods of enumerative geometry. I. Vainsencher \cite {V} tried
to use various blow-ups of the ordinary $ \PP^5 $ of conics as moduli spaces,
but the intersection of the five tangency conditions in this moduli space
always resulted in a scheme with many non-enumerative components, most of which
had too big dimension, were non-reduced or even had embedded components. Their
geometry was so complicated that the problem could not be solved that way.

In this paper we use different moduli spaces, namely moduli spaces of relative
stable maps, to solve the problem. There is a well-defined compact moduli space
$ \bar M_{(2,2,2,2,2)}^Y (\PP^2,2) \subset \bar M_{0,5} (\PP^2,2) $ that
parametrizes rational stable maps to $ \PP^2 $ of degree 2 (i.e.\ conics) with
5 marked points such that the stable map is tangent to $Y$ at all these
points. It comes equipped with a 0-dimensional virtual fundamental class, whose
degree $ N_d $ can be computed explicitly using the methods of \cite {G}.

We can interpret the number $ N_d $ as the ``virtual number'' of conics that
are 5-fold tangent to $Y$. It is only virtual because it contains --- just as
in Vainsencher's classical computations --- non-enumerative contributions from
the ``boundary'' of the moduli space. These contributions are quite simple
however. It is not hard to see that the only degree-2 rational stable maps $ f:
C \to \PP^2 $ that satisfy the tangency conditions at the 5 marked points are
all double covers of a bitangent of $Y$, and have the marked points distributed
in one of the following two ways:
\begin {center} \epsfig {file=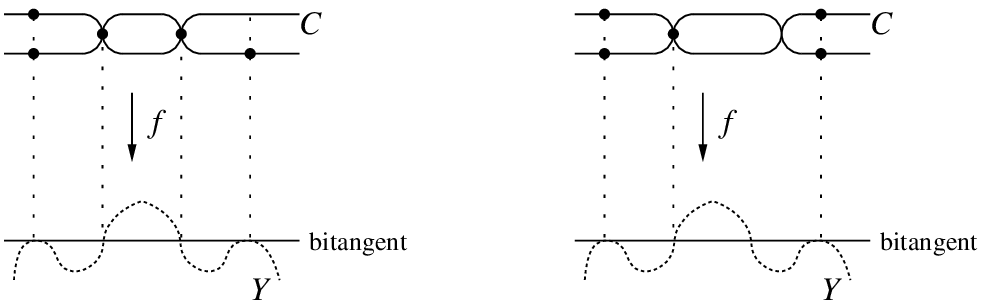} \end {center}
There are only finitely many stable maps with marked points as in this picture
on the left, so we just have to count them and subtract their number from the
virtual invariant $ N_d $. The picture on the right however shows a 1-dimensional
family of stable maps (the unmarked ramification point of $f$ can move). We
will use equations from relative Gromov-Witten theory to compute the degree of
the 0-dimensional virtual fundamental class of the moduli space on this
1-dimensional component. By subtracting both correction terms, we finally
arrive at the enumerative numbers $ n_d $. They are
    \[ n_d = \frac 1 {5!} \, d \, (d-3) \, (d-4) \,
             (d^7+12d^6-18d^5-540d^4+251d^3+5712d^2-1458d-14580). \]
Our (as well as Vainsencher's) motivation for studying this problem came from a
question concerning rational curves in K3 surfaces. If $X$ is a K3 surface and
$ \beta \in H_2 (X,\ZZ) $ the homology class of a holomorphic curve in $X$, then
various authors \cite {YZ}, \cite {B}, \cite {Go}, \cite {BL} have shown that
the number of rational curves in $X$ of class $ \beta $ is equal to the $ q^d $
coefficient of the series
  \[ G(q) = \prod_{i>0} \frac 1 {(1-q^i)^{24}}
          = 1 + 24q + 324q^2 + 3200q^3
              + 25650q^4 + 176256q^5 + \cdots \]
with $ d = \frac 1 2 \beta^2 + 1 $, \emph {if the class $ \beta $ is
primitive}, i.e.\ not a non-trivial multiple of a smaller integral homology
class. There is a well-defined ``K3 invariant'' (using a modified
obstruction theory on the moduli spaces of stable maps to $X$) for
non-primitive $ \beta $ too \cite {BL}; it is however not known yet how this
invariant relates to the above series $ G(q) $ or to the enumerative number.

The result of this paper allows us to study this question in a non-trivial
example: we let $X$ be the double cover of $ \PP^2 $ branched along a sextic
curve $Y$, and take $ \beta $ to be the pull-back of the class of conics in $
\PP^2 $. Our work allows us to compute the enumerative number of rational
curves in $X$ of class $ \beta $, which we can then compare to the
corresponding number of the series $ G(q) $. The result is that the ``K3
invariant'' is equal to the corresponding term in the $ G(q) $ series, plus a
double cover correction term that is equal to $ \frac 1 8 $ times the number of
rational curves in $X$ of class $ \frac 1 2 \beta $. Note that this is the same
sort of correction term as for multiple covers of rational curves in Calabi-Yau
threefolds. We conjecture that this pattern continues for classes $ \beta $ of
higher divisibility.

The paper is organized as follows. In section \ref {sec-relative} we will show
how to compute the relative Gromov-Witten invariant $ N_d $. We analyze the
moduli space $ \bar M_{(2,2,2,2,2)}^Y (\PP^2,2) $ and its virtual fundamental
class in sections \ref {sec-components} and \ref {sec-virtual}, respectively,
leading to the final result for $ n_d $ in corollary \ref {nd}. Section \ref
{sec-k3} contains the application to K3 surfaces mentioned above.

We thank J. Kock for introducing us to the problem, as well as for a Maple
program to compute the relative Gromov-Witten invariants using the algorithm of
\cite {G}. The author is also grateful to the Institute for Advanced Study for
its hospitality and stimulating working atmosphere. This work has been funded
by the NSF grant DMS 9729992.


\section {Relative Gromov-Witten invariants} \label {sec-relative}

In this section we will show how to compute the relative Gromov-Witten
invariant that corresponds to the number of conics that are 5-fold tangent to a
given smooth plane curve. We will use the notations and results from \cite {G},
to which we also refer for further details.

\begin {definition}
  Let $ Y \subset \PP^2 $ be a smooth curve of degree $d$, and let $
  m_1,\dots,m_n $ be non-negative integers. We denote by $ \bar
  M_{(m_1,\dots,m_n)} = \bar M^Y_{(m_1, \dots,m_n)} (\PP^2,2) $ the moduli
  space of $n$-pointed relative stable maps of degree 2 to $ \PP^2 $ relative
  to $Y$ with multiplicities $ m_1,\dots,m_n $, as defined in \cite {G}
  definitions 1.1 and 1.18. Its virtual fundamental class is denoted $ [\bar
  M_{(m_1,\dots,m_n)}]^{virt} $.
\end {definition}

\begin {remark}
  The moduli space $ \bar M_{(m_1,\dots,m_n)} $ can be thought of as a
  compactification of the space of irreducible plane conics together with $n$
  distinct marked points on them at which the conic has the prescribed local
  intersection multiplicities $ m_1,\dots,m_n $ with $Y$. In particular, the
  moduli space $ \bar M_{(2,2,2,2,2)} $ corresponds to conics 5-fold tangent to
  $Y$.
\end {remark}

\begin {remark} \label {rel-def}
  For future reference let us recall the precise definition from \cite
  {G}. Consider the degree-$d$ Veronese embedding $ i: \PP^2 \to \PP^D $ with $
  D = \binom {D+2}{2}-1 $. We have $ i(Y) = i(\PP^2) \cap H $ in $ \PP^D $ for
  a suitable hyperplane $ H \subset \PP^D $. The inclusion morphism $i$ induces
  an inclusion of moduli spaces $ \bar M_{0,n} (\PP^2,2) \subset \bar M_{0,n}
  (\PP^D,2d) $. Moreover, let $ \bar M^H_{(m_1,\dots,m_n)} (\PP^D,2d) $ be the
  closure in $ \bar M_{0,n} (\PP^D,2d) $ of all stable maps $
  (C,x_1,\dots,x_n,f) $ such that $C$ is irreducible, $ f(C) \not\subset H $,
  and the divisor $ f^* Y $ on $C$ contains the points $ x_i $ with
  multiplicities $ m_i $. Then we define
    \[ \bar M_{(m_1,\dots,m_n)} := \bar M_{0,n} (\PP^2,2) \cap
         \bar M^H_{(m_1,\dots,m_n)} (\PP^D,2d) \]
  (with the intersection taken in $ \bar M_{0,n} (\PP^D,2d) $). The virtual
  fundamental class is the corresponding refined intersection product
    \[ [\bar M_{(m_1,\dots,m_n)}]^{virt} := [\bar M_{0,n} (\PP^2,2)] \cdot
         [\bar M^H_{(m_1,\dots,m_n)} (\PP^D,2d)]
         \in A_* (\bar M_{(m_1,\dots,m_n)}). \]
  The virtual dimension of $ \bar M_{(m_1,\dots,m_n)} $ is $ 5 - \sum_i (m_i-1)
  $. By abuse of notation, we will always drop the superscript $ virt $ from
  the notation of the virtual fundamental class from now on, as we do not need
  the ordinary fundamental classes of these spaces.
\end {remark}

\begin {remark} \label {easy-mult}
  From remark \ref {rel-def} we get immediately the following statement: let $
  \calC = (C,x_1,\dots,x_n,f) \in \bar M_{(m_1,\dots,m_n)} $ be an
  automorphism-free stable map such that $ C \cong \PP^1 $ is irreducible and $
  f(C) \not\subset Y $. Then, locally around this point, $ \bar
  M_{(m_1,\dots,m_n)} $ is scheme-theoretically the subscheme of $ \bar M_{0,n}
  (\PP^2,2) $ given by the $ \sum_i m_i $ equations that describe the vanishing
  of the $ m_i $-jets of $ ev_i^* Y $ at the points $ x_i $, where $ ev_i: \bar
  M_{0,n} (\PP^2,2) \to \PP^2 $ are the evaluation maps. If moreover $ \bar
  M_{(m_1,\dots,m_n)} $ has the expected dimension at this point $ \calC $,
  then this point lies on a unique irreducible component of $ \bar
  M_{(m_1,\dots,m_n)} $, and the virtual fundamental class on this component is
  just the ordinary scheme-theoretic fundamental class, i.e.\ the length of the
  scheme $ \bar M_{(m_1,\dots,m_n)} $ at this irreducible component.
\end {remark}

\begin {remark} \label {set-descr}
  There is an easier description of $ \bar M_{(m_1,\dots,m_n)} $ \emph {as a
  set}. Namely, $ \bar M_{(m_1,\dots,m_n)} $ is the subspace of $ \bar M_{0,n}
  (\PP^2,2) $ of all $n$-pointed rational stable maps $ (C,x_1,\dots,x_n,f) $
  of degree 2 to $ \PP^2 $ such that the following two conditions are
  satisfied:
  \begin {enumerate}
  \item $ f(x_i) \in Y $ for all $i$ such that $ m_i > 0 $,
  \item $ f^* Y - \sum_i m_i x_i \in A_0 (f^{-1}(Y)) $ is effective.
  \end {enumerate}
\end {remark}

As mentioned above, the moduli space $ \bar M_{(2,2,2,2,2)} $ has virtual
dimension zero and corresponds to conics 5-fold tangent to $Y$ (together with a
labeling of the five tangency points). Hence we define:

\begin {definition}
  The number
    \[ N_d := \frac 1 {5!} \cdot \deg [\bar M_{(2,2,2,2,2)}] \in \QQ \]
  will be called the virtual number of conics 5-fold tangent to $Y$.
\end {definition}

The number $ N_d $ is only virtual because it receives correction terms from
double covers of lines (see sections \ref {sec-components} and \ref
{sec-virtual}). In the rest of this section we will show how to compute the
number $ N_d $. Obviously, we can assume that $ d \ge 5 $.

The computation is done using the main theorem 2.6 of \cite {G} that tells us
``how to raise the multiplicities of the moduli spaces'': it says that
\begin {equation} \label {main-eqn}
  (ev_n^* Y + m_n \psi_n) \cdot [\bar M_{(m_1,\dots,m_n)}] =
     [\bar M_{(m_1,\dots,m_{n-1},m_n+1)}] +
       \mbox {correction terms},
\end {equation}
where $ ev_n: \bar M_{(m_1,\dots,m_n)} \to \PP^2 $ is the evaluation map at the
last marked point, and $ \psi_n $ is the first Chern class of the cotangent
line bundle $ L_n $, i.e.\ of the bundle whose fiber at a stable map $
(C,x_1,\dots,x_n,f) $ is the cotangent space $ T^\vee_{C,x_n} $. The correction
terms are as follows. Every correction term corresponds to a moduli space of
reducible curves with $ r+1 $ components $ C_0,\dots,C_r $, where $ C_0 $ is
contracted\footnote {This uses the fact that the curve $Y$ has positive genus,
and that therefore every rational stable map to $Y$ must be constant. In
general, $ C_0 $ can be a curve with any homology class in $Y$.} to a point of
$Y$, and the other components $ C_i $ intersect $ C_0 $ in a point where they
have local intersection multiplicity $ \mu_i $ to $Y$. We get such a correction
term for every $r$, every choice of the $ \mu_i $, and every splitting of the
total homology class and the marked points onto the components $ C_i $, such
that the following two conditions are satisfied:
\begin {enumerate}
\item [(a)] the last marked point $ x_n $ lies on the component $ C_0 $,
\item [(b)] the sum of all the $ \mu_i $ is equal to the sum of those $ m_i $
  such that $ x_i \in C_0 $.
\end {enumerate}
These correction terms appear in the above equation with multiplicity $
\prod_{i=1}^r \mu_i $.

\begin {example} \label {last-eqn}
  Here is an example of equation (\ref {main-eqn}). In the case
    \[ (ev_5^* Y + \psi_5) \cdot [\bar M_{(2,2,2,2,1)}] =
         [\bar M_{(2,2,2,2,2)}] +
         \mbox {correction terms} \]
  we want to figure out the correction terms. As this is an equation in
  (virtual) dimension 0, the contracted component $ C_0 $ must have exactly 3
  special points (it would not be stable if it had fewer, and it would have
  moduli if it had more). Hence the correction terms fall into these two
  categories:
  \begin {enumerate}
  \item $ r=1 $ (in the above notation), $ C_1 $ is a conic, and $ C_0 $ is a
    contracted component with three special points $ x_5 $, the intersection
    point with $ C_1 $, and one other $ x_i $ for $ i=1,\dots,4 $,
  \item $ r=2 $, $ C_1 $ and $ C_2 $ two lines, and $ C_0 $ is a contracted
    component with three special points $ x_5 $ and the two intersection points
    with $ C_1 $ and $ C_2 $.
  \end {enumerate}
  Actually, case (ii) cannot occur, because condition (b) above cannot be
  satisfied: the sum $ \mu_1 + \mu_2 $ is at least 2, whereas $ m_5 $ is only
  1. Hence the only correction terms are of type (i). We get four of them: one
  for each choice of the point $ x_i $ that is to lie on the contracted
  component $ C_0 $. We have $ \mu_1 = m_i + m_5 = 3 $ in each of these cases
  by condition (b). All four correction terms appear with multiplicity $ \mu_1
  = 3 $. Pictorially, the equation reads
  \begin {center} \epsfig {file=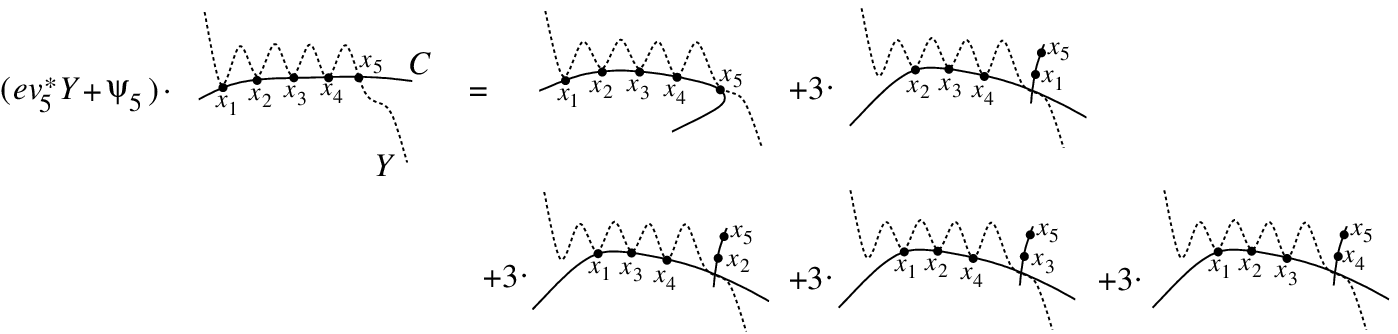,width=125mm} \end {center}
  Here, the dotted curve is the fixed curve $Y$, and the solid curve is the
  moving conic $C$. In the four correction terms, the component with $ x_5 $
  on it is meant to be contracted. Written down as an equation of virtual
  fundamental classes of moduli spaces, the formula reads
  \begin {align} \label {last}
    (ev_5^* Y + \psi_5) \cdot [\bar M_{(2,2,2,2,1)}] = &
      [\bar M_{(2,2,2,2,2)}] +
      3 [\bar M_{(3,2,2,2)}] \\
    & + 3 [\bar M_{(2,3,2,2)}] + 3 [\bar M_{(2,2,3,2)}]
      + 3 [\bar M_{(2,2,2,3)}]. \notag
  \end {align}
\end {example}

\begin {remark}
  The correction terms in equation (\ref {main-eqn}) are themselves products of
  moduli spaces of relative stable maps of the form $ \bar
  M_{(m_1',\dots,m_k')}^Y (\PP^2,d') $ with $ d' \le 2 $ and $ m_1' + \cdots
  m_k' \le m_1 + \cdots m_n $. In other words, this equation expresses
  invariants (i.e.\ intersection products of $ ev_i^* H $ and $ \psi_i $
  classes, where $H$ is a line in $ \PP^2 $) on the relative moduli space $
  \bar M_{(m_1,\dots,m_n+1)} $ in terms of other invariants on relative moduli
  spaces $ \bar M_{(m_1',\dots,m_k')}^Y (\PP^2,d') $ whose ``total
  multiplicity'' $ \sum_i m_i' $ is smaller than the total multiplicity $
  1+\sum_i m_i $ of $ M_{(m_1,\dots,m_n+1)} $. Hence, applying equation (\ref
  {main-eqn}) recursively $ m_1 + \cdots + m_n $ times, we can express every
  invariant on $ \bar M_{(m_1,\dots,m_n)} $ in terms of invariants on $ \bar
  M_{(0,\dots,0)} $, which are just ordinary Gromov-Witten invariants of $
  \PP^2 $. As the Gromov-Witten invariants of $ \PP^2 $ are well-known, we can
  thus compute all relative Gromov-Witten invariants recursively, in particular
  $ N_d $. Example \ref {last-eqn} is the first step in this recursion process;
  it expresses the invariant $ \bar M_{(2,2,2,2,2)} $ (with total multiplicity
  10) in terms of invariants with total multiplicity 9.
\end {remark}

Without actually carrying out the recursion, we can see the following.

\begin {lemma}
  The function $ d \mapsto N_d $ is a polynomial of degree 10 with leading
  coefficient $ \frac 1 {5!} $.
\end {lemma}

\begin {proof}
  Using equation (\ref {main-eqn}) it is easy to show by induction that every
  invariant (i.e.\ intersection product of $ ev_i^* H $ and $ \psi_i $ classes)
  on a moduli space $ \bar M_{(m_1,\dots,m_n)} $ is a polynomial in $d$ of
  degree (at most) $ m_1 + \cdots + m_n $. In fact, this is obvious for $ m_1 +
  \cdots + m_n = 0 $, as we then just have ordinary Gromov-Witten invariants of
  $ \PP^2 $ (that do not depend on $Y$). Equation (\ref {main-eqn}) reads
    \[ [\bar M_{(m_1,\dots,m_{n-1},m_n+1)}] =
       (d \; ev_n^* H + m_n \psi_n) \cdot [\bar M_{(m_1,\dots,m_n)}]
         - \mbox {correction terms}. \]
  All correction terms have total multiplicity at most $ m_1 + \cdots + m_n $,
  so by induction hypothesis they contribute a polynomial in $d$ of degree
  at most $ m_1 + \cdots + m_n $. The same is true for the $ \psi_n $ summand
  on the right hand side. Hence, as every invariant on $ \bar M_{(m_1,\dots,
  m_n)} $ is a polynomial in $d$ of degree at most $ m_1 + \cdots + m_n $ by
  assumption, it follows that every invariant on $ \bar M_{(m_1,\dots,m_n+1)}
  $ is a polynomial in $d$ of degree at most $ m_1 + \cdots + m_n + 1 $.

  It can be seen from the same recursive formula that the $ d^{10} $
  coefficient of the invariant $ \deg [\bar M_{(2,2,2,2,2)}] $ is just
    \[ \prod_{i=1}^5 ev_i^* H^2 \cdot [\bar M_{0,5} (\PP^2,2)], \]
  i.e.\ the number of conics through 5 general points in the plane. This number
  is 1, proving the statement of the lemma about the leading coefficient.
\end {proof}

The precise form of the polynomial $ N_d $ is quite complicated and can only be
obtained by carrying out the full recursion as described above. We only give
the result here; a Maple program to compute it can be obtained from the author
on request.

\begin {proposition} \label {Nd}
  For $ d \ge 5 $, the virtual number of conics 5-fold tangent to $Y$ is
    \[ N_d = \frac 1 {5!} \, d \, (d-3) \, (d-4) \,
             (d^7+12d^6-18d^5-540d^4+311d^3+5457d^2-2133d-12690). \]
\end {proposition}

\begin {remark}
  The first few values of $ N_d $ are given in the following table.
    \[ \begin {array}{|c|c|c|c|c|c|c|} \hline
         d & 5 & 6 & 7 & 8 & 9 & 10 \\ \hline
         N_d & 1985 & 71442 & 687897 & 3893256 & 16180398 & 54679380 \\ \hline
       \end {array} \]
\end {remark}

\begin {remark} \label {refined}
  In this section we have only used equation (\ref {main-eqn}) in the Chow ring
  of the moduli space of stable maps $ \bar M_{0,n} (\PP^2,2) $. In fact, there
  is a refined version of this equation that we will need in section \ref
  {sec-virtual}. If $ \calP_k $ denotes the rank-$ (k+1) $ bundle of (relative)
  $k$-jets of $ ev_n^* \calO(Y) $, there is a section $ \sigma $ of the line
  bundle $ \calP_{m_n} / \calP_{m_n-1} $ on $ \bar M_{(m_1,\dots,m_n)} $ whose
  vanishing describes precisely the condition that the map $f$ of a stable map
  $ (C,x_1,\dots,x_n,f) $ has multiplicity (at least) $ m_n+1 $ to $Y$ at $ x_n
  $. The first Chern class of this line bundle is $ ev_n^* Y + m_n \psi_n
  $. (In fact, this is the idea how equation (\ref {main-eqn}) is proven.) The
  refined version of equation (\ref {main-eqn}) now states that this equation
  also holds in the Chow group of the zero locus of the section $ \sigma $ on $
  \bar M_{(m_1,\dots,m_n)} $.
\end {remark}


\section {The components of $ \bar M_{(2,2,2,2,2)} $}
  \label {sec-components}

Having just computed the invariant $ N_d $, we will now study its enumerative
significance. To do this, we have to identify the components of $ \bar
M_{(2,2,2,2,2)} $ and compute their virtual fundamental classes. We will
assume from now on that $Y$ is generic of degree $ d \ge 5 $. The equation of
$Y$ is $ F = \sum_I a_I z^I = 0 $, where $I$ runs over all multi-indices $
(i_0,i_1,i_2) $ with $ i_0+i_1+i_2=d $, and $ z_0,z_1,z_2 $ are the homogeneous
coordinates on $ \PP^2 $. We start with irreducible stable maps whose image is
a smooth conic.

\begin {lemma} \label {expdim-irr}
  Let $ m_1,\dots,m_n $ be non-negative integers such that $ \sum m_i \le 10
  $, and let $ \calC \in \bar M_{(m_1,\dots,m_n)} $ be an irreducible stable map
  that is not a double cover of a line. Then:
  \begin {enumerate}
  \item The moduli space $ \bar M_{(m_1,\dots,m_n)} $ is smooth of dimension $
    5 - \sum_i (m_i-1) $ at $ \calC $ (which is the expected dimension).
  \item If this expected dimension is negative, then there is no such point $
    \calC $.
  \end {enumerate}
\end {lemma}

\begin {proof}
  The plane degree-$d$ curves are parametrized by a projective space $ \PP^D $
  with $ D = \frac 1 2 d(d+3) $, whose coordinates are the coefficients $ a_I $
  of $F$. Let $ Z \subset \bar M_{0,n}(\PP^2,2) \times \PP^D $ be the closed
  substack of pairs $ ((C,x_1,\dots,x_n,f),Y) $ such that the pull-back by $f$
  of the equation of $Y$ vanishes at the points $ x_i $ to order $ m_i $ for
  all $i$. We claim that $Z$ is smooth of the expected dimension at every point
  $ ((C,x_1,\dots,x_n,f),Y) $ such that $C$ is irreducible and $f$ is not a
  double cover of a line.

  To prove this, we have to show that the matrix of derivatives of the
  equations describing $Z$ has maximal rank at the given point $
  ((C,x_1,\dots,x_n,f),Y) $. By a projective coordinate transformation of $
  \PP^2 $ and choosing homogeneous coordinates on $ C \cong \PP^1 $, we can
  assume that the map $f$ is given by $ (s:t) \mapsto (s^2:st:t^2) $, and the
  $n$ marked points are $ (1:\lambda_i) $ with pairwise distinct $ \lambda_i $.

  Let us now write down the derivatives of the multiplicity equations with
  respect to the first $ m := \sum_i m_i $ of the variables $ a_{(d,0,0)} $, $
  a_{(d-1,1,0)} $, $ a_{(d-1,0,1)} $, $ a_{(d-2,1,1)} $, $ a_{(d-2,0,2)} $, $
  a_{(d-3,1,2)} $, $ a_{(d-3,0,3)} $, $ a_{(d-4,1,3)} $, $ a_{(d-4,0,4)} $, $
  a_{(d-5,1,4)} $ (remember $ m \le 10 $ and $ d \ge 5 $). These coordinates
  are chosen to be the coefficients of $ s^{2d-i} t^i $ for $ i=0,\dots,9 $
  when we substitute the map $f$ into $F$.

  Multiplicity $ m_i $ at the point $ (s:t)=(1:\lambda_i) $ means that $
  F|_{s=1,t=\lambda_i+\epsilon} $ has no $ \epsilon $ terms of order less than $
  m_i $. So the rows of the derivative matrix are just $ (\binom {i}{k}
  \lambda_j^{i-k})_{i=0,\dots,m-1} $, for $ 0 \le k < m_j $ and $ 1 \le j \le n
  $. For example, for $ m_1 = m_2\ = m_3 = m_4 = m_5 = 2 $ we get the matrix
    \[ \left( \begin {array}{ccccc}
         1 & \lambda_1 & \lambda_1^2 & \cdots & \lambda_1^9 \\
         0 & 1 & 2\lambda_1 & \cdots & 9\lambda_1^8 \\
         \vdots & \vdots & \vdots & \vdots & \vdots \\
         1 & \lambda_5 & \lambda_5^2 & \cdots & \lambda_5^9 \\
         0 & 1 & 2\lambda_5 & \cdots & 9\lambda_5^8 \\
       \end {array} \right). \]
  By subtracting $ \lambda_1 $ times the $i$-th column from the $ (i+1) $-st
  column for $ 1 \le i < m $ and using induction, we see that the determinant
  is $ \prod_{i<j} (\lambda_i - \lambda_j)^{m_i m_j} $. In particular, it is
  not zero, so $Z$ is smooth of the expected dimension at $
  ((C,x_1,\dots,x_n,f),Y) $.

  By remark \ref {easy-mult}, the statement of the lemma is just that the fiber
  of $Z$ over a general point of $ \PP^D $ is smooth of the expected dimension
  around a point considered above. This follows now from the Bertini theorem.
\end {proof}

Using remark \ref {easy-mult} again, the following two corollaries are immediate.

\begin {corollary}
  Every irreducible stable map in $ \bar M_{(m_1,\dots,m_n)} $ whose image in $
  \PP^2 $ is a smooth conic lies in a unique irreducible component of $ \bar
  M_{(m_1,\dots,m_n)} $ of the expected dimension. The virtual fundamental
  class of this component is equal to the usual one.
\end {corollary}

\begin {corollary}
  The number of smooth plane conics 5-fold tangent to $Y$ is finite. We denote
  it by $ n_d $.
\end {corollary}

We will now study the additional non-enumerative contributions to the virtual
invariant $ N_d $.

\begin {lemma}
  The moduli space $ \bar M_{(2,2,2,2,2)} $ has the following connected
  components:
  \begin {enumerate}
  \item [(A)] $ 5! $ points for every smooth conic 5-fold tangent to $Y$.
  \item [(B)] $ 5! \cdot (d-4)(d-5) $ points for every bitangent of $Y$,
    corresponding to double covers of the bitangent, with marked points as in
    the picture below on the left, i.e.\ the map is ramified over two
    transverse intersection points of the bitangent with $Y$, and the five
    marked points are the two ramification points, both inverse image points
    of one bitangency point, and one inverse image point of the other
    bitangency point.
    \begin {center} \epsfig {file=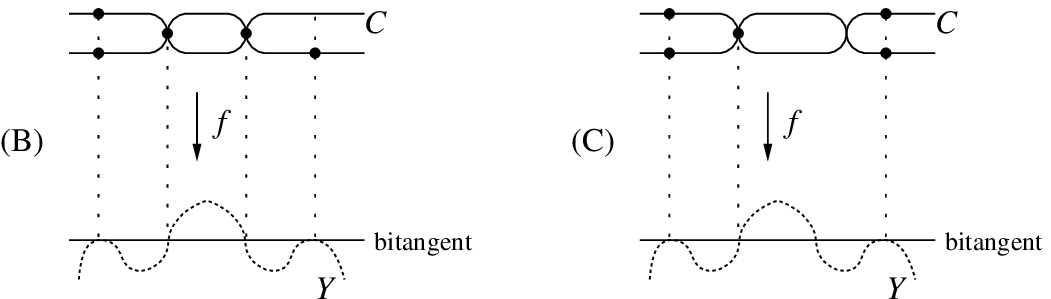} \end {center}
  \item [(C)] $ \frac 1 4 \cdot 5! \cdot (d-4) $ smooth rational curves for
    every bitangent of $Y$, corresponding to double covers of the bitangent,
    with marked points as in the picture above on the right, i.e.\ (for a
    general stable map in this smooth rational curve) the map is ramified over
    one transverse intersection point of the bitangent with $Y$ and one other
    arbitrary point, and the five marked points are the first ramification
    point, and the four inverse image points of the two bitangency points.
  \end {enumerate}
\end {lemma}

\begin {proof}
  Case 1: the image of the stable map is a smooth conic. If the five marked
  points are distinct in $ \PP^2 $, we get the components (A) by lemma \ref
  {expdim-irr}, with the $ 5! $ corresponding to the labeling of the marked
  points. If two of the points coincide in $ \PP^2 $ (i.e.\ lie on a contracted
  component of the stable map), then by the description of $ \bar
  M_{(2,2,2,2,2)} $ in remark \ref {set-descr} the conic must have contact of
  order (at least) 4 to $Y$ at this point, i.e.\ it lies in $ \bar
  M_{(4,2,2,2)} $. But this space is empty by lemma \ref {expdim-irr}.

  Case 2: the image of the stable map is a union of two (distinct) lines. It is
  easy to see that the conditions of remark \ref {set-descr} cannot be
  satisfied in this case.

  Case 3: the stable map is a double cover of a line. There are six possible
  points of tangency to $Y$: the four inverse image points of the bitangency
  points, and the two ramification points if they are mapped to points of
  $Y$. For the stable map in $ \bar M_{(2,2,2,2,2)} $ we can pick any five of
  these six points. If we leave out one of the points over the bitangency
  points, we arrive at the components (B), otherwise we get the components
  (C).

  In case (B) we get a factor of $ 5! $ for the choice of labeling of the
  marked points, a factor of $ \binom {d-4}{2} $ for the choice of two
  transverse intersection points of $Y$ with the bitangent, and a factor of 2
  for the choice of bitangency point over which we take only one inverse image
  point to be marked. In case (C) the second ramification point is not fixed,
  so we get one-dimensional families of such curves. Every such family has a
  2:1 map to the bitangent given by the image of the moving ramification point;
  the two stable maps in a fiber of this map differ by exchanging the marked
  points over one bitangency point. The map is simply ramified over the two
  stable maps where the moving ramification point is one of the bitangency
  points. Hence every such family is a $ \PP^1 $. The number of such families
  is $ d-4 $ (for the choice of transverse intersection point of the bitangent
  with $Y$) times $ \frac 1 4 \cdot 5! $ (for the labeling of the marked
  points, taking into account that exchanging the marked points over the
  bitangency points does not give us a new family).
\end {proof}

Of course, the virtual fundamental class of $ \bar M_{(2,2,2,2,2)} $ splits
naturally into a sum of virtual fundamental classes on each of the connected
components that we have just identified. As it is well-known that the number of
bitangents of $Y$ is $ \frac 1 2 d (d+3)(d-2)(d-3) $ (see e.g.\ \cite {H}
exercise IV.2.3 f), we get the following corollary.

\begin {corollary} \label {nd-cor}
  We have
    \[ N_d = n_d + \frac 1 2 d (d+3)(d-2)(d-3)(d-4)(d-5) b_d
                 + \frac 1 8 d (d+3)(d-2)(d-3)(d-4) c_d, \]
  where $ b_d $ is the degree of the part of the virtual fundamental class of
  $ \bar M_{(2,2,2,2,2)} $ supported on the point $ \bar M_B \in \bar
  M_{(2,2,2,2,2)} $ below, and $ c_d $ is the corresponding degree supported on
  the smooth rational curve $ \bar M_C \subset \bar M_{(2,2,2,2,2)} $ below.
  \begin {center} \epsfig {file=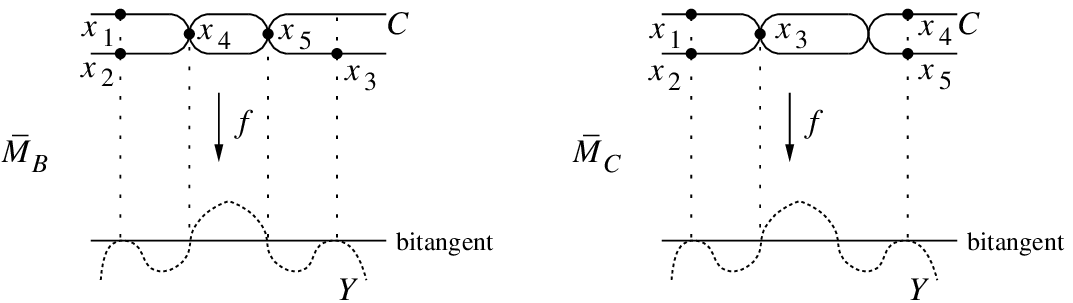} \end {center}
\end {corollary}


\section {Computation of the virtual fundamental classes}
  \label {sec-virtual}

In this section we will do the necessary computations to determine the numbers
$ b_d $ and $ c_d $ of corollary \ref {nd-cor}. Most of them are simple
calculations in local coordinates, so we will only sketch these parts and leave
the details to the reader.

The computation of $ b_d $ is quite simple, as the component $ \bar M_B $ of $
\bar M_{(2,2,2,2,2)} $ has the expected dimension.

\begin {lemma} \label {bd}
  $ b_d = 1 $ for all $d$.
\end {lemma}

\begin {proof}
  By remark \ref {easy-mult} we just have to show that the 10 equations of
  vanishing of the 1-jets of$ ev_i^* F $ for $ i=1,\dots,5 $ locally cut out
  the point $ \bar M_B $ in $ \bar M_{0,5} (\PP^2,2) $ scheme-theoretically
  with multiplicity 1. Let us start with the 1-jets at the points $ x_1,x_2,x_3
  $, i.e.\ with the space $ \bar M_{(2,2,2)} $. We can choose the coordinates on $
  \PP^2 $ such that the bitangent is $ \{ z_2=0 \} \subset \PP^2 $ and the
  bitangency points are $ (1:0:0) $ and $ (0:1:0) $. This means that
  \begin {align*}
    a_{(d,0,0)} = a_{(d-1,1,0)} = 0, &
     \quad a_{(d-1,0,1)}, a_{(d-2,2,0)} \neq 0
     \qquad \mbox {(tangency at $ (1:0:0) $),} \\
    a_{(0,d,0)} = a_{(1,d-1,0)} = 0, &
     \quad a_{(0,d-1,1)}, a_{(2,d-2,0)} \neq 0
     \qquad \mbox {(tangency at $ (0:1:0) $).}
  \end {align*}
  Moreover, we can choose coordinates on the source $ \PP^1 $ such that the
  stable map $ \bar M_B $ is given by $ (s:t) \mapsto (s^2-t^2:st:0) $, and the
  marked points are $ x_1 = (1:0) $, $ x_2 = (0:1) $, $ x_3 = (1:1) $. Local
  coordinates of $ \bar M_{0,3} (\PP^2,2) $ around this point are then $
  \epsilon_1,\dots,\epsilon_8 $, where the stable map is given by
    \[ (s:t) \mapsto (s^2-t^2+\epsilon_1 s^2 + \epsilon_2 st + \epsilon_3 t^2:
                      st+\epsilon_4 s^2 + \epsilon_5 t^2:
                      \epsilon_6 s^2 + \epsilon_7 st + \epsilon_8 t^2). \]
  The three tangency equations are that $ F|_{s=1,t=\xi} $, $ F|_{s=\xi,t=1}
  $, and $ F|_{s=1,t=1+\xi} $ have no constant and linear $ \xi $ terms. It
  is an easy computation to see that these 6 equations, linearized in the $
  \epsilon_i $, give
    \[ \epsilon_1 + \epsilon_2 + \epsilon_3 = \epsilon_4 = \epsilon_5 =
       \epsilon_6 = \epsilon_7 = \epsilon_8 = 0,\]
  so $ \bar M_{(2,2,2)} $ is smooth of dimension 2 at the point $ \bar M_B $
  (with the points $ x_4 $ and $ x_5 $ forgotten).

  Now let us consider the two other tangency conditions at the points $ x_4 $
  and $ x_5 $. As the four coordinates of $ \bar M_{(2,2,2,0,0)} $ around $
  \bar M_B $ we can choose the images of the ramification points and a point in
  the domain of the stable map in the neighborhood of each ramification
  point. Considering only one ramification point for now, the two corresponding
  local coordinates are $ \tilde \epsilon_1 $ and $ \tilde \epsilon_2 $, where
  the stable map is given locally in affine coordinates as $ t \mapsto t^2 +
  \tilde \epsilon_1 $, and the marked point is $ t = \tilde \epsilon_2
  $. Tangency means that the constant and linear $ \xi $ terms of $ (\tilde
  \epsilon_2 + \xi)^2 + \tilde \epsilon_1 $ vanish, so linearly in the $ \tilde
  \epsilon_i $ we get $ \tilde \epsilon_1 = \tilde \epsilon_2 = 0 $. The same
  is true for the other ramification point, so we see that $ \bar M_B $ is a
  smooth point of $ \bar M_{(2,2,2,2,2)} $.
\end {proof}

To study the space $ \bar M_C $, we need a lemma that tells us how the stable
maps in $ \bar M_C $ can be deformed if we relax some of the multiplicity
conditions.

\begin {lemma} \label {no-def}
  Let $H$ be a line in $ \PP^2 $, and let $ P \in H $ be a point where $Y$ is
  simply tangent to $H$. Let $ \calC = (C,x_1,x_2,f) \in \bar M_{(2,2)} $ be a
  (possibly reducible) double cover of $H$, such that $ f^{-1} (P) = \{ x_1,x_2
  \} $. Then every stable map in $ \bar M_{(2,2)} $ in a neighborhood of $
  \calC $ is also a double cover of a (maybe different) line.
\end {lemma}

\begin {proof}
  It is obvious that $ \calC $ cannot be deformed into a union of two distinct
  lines in $ \bar M_{(2,2)} $. So we have to show that $ \calC $ cannot be
  deformed to an irreducible smooth conic in $ \bar M_{(2,2)} $.

  We use the classical space of complete conics (which is isomorphic to $ \bar
  M_{0,0} (\PP^2,2) $). Recall that this space is the closure in $ \PP^5 \times
  (\PP^5)^\vee $ of the set $ (C,C^\vee) $, where $C$ is an irreducible conic
  and $ C^\vee $ its dual. For our given point $ \calC $ (with the two marked
  points forgotten for the moment), $C$ is the double line $H$, and $ C^\vee $
  is the union of the two lines in $ \PP^\vee $ that correspond to the two
  ramification points of $f$ in $ \PP^2 $. Assume that we can deform $
  (C,C^\vee) $ in the space of complete conics to an irreducible conic that is
  still tangent to $Y$ at two points in the neighborhood of $P$. In particular,
  we would then deform $ C^\vee $ to an irreducible conic that is tangent to
  the dual of $Y$ (or more precisely if $H$ is a bitangent of $Y$: to the
  branch of the dual of $Y$ that corresponds to the point $P$) at two points in
  this neighborhood. By the continuity of intersection products this means that
  both lines of $ C^\vee $ must actually be the line corresponding to the point
  $P$. Hence both ramification points of $f$ would have to be $P$. This means
  that $ \calC $ must have a contracted rational component over $P$, in
  contradiction to the assumption $ f^{-1} (P) = \{ x_1,x_2 \} $.
\end {proof}

We want to reduce the computation of $ c_d $ to spaces that have the expected
dimension. To do this, we use equation (\ref {last}) from example \ref
{last-eqn}. Note that by remark \ref {refined} this equation is true in the
Chow group of the geometric zero locus of the section $ \sigma $ (that
describes the tangency condition and whose zero locus has class $ ev_5^* Y +
\psi_5 $), so it makes sense to restrict the equation to a connected component
of this zero locus. Using lemma \ref {no-def} it is easy to see that $ \bar M_C
$ is such a connected component, so we will restrict equation (\ref {last}) to
$ \bar M_C $ and denote this restriction by $ |_{\bar M_C} $.

Note first that on $ \bar M_C $ the point $ x_5 $ can only come close to $ x_4
$, but never to the other three marked points. Hence equation (\ref {last})
restricted to $ \bar M_C $ reads
\begin {equation} \label {last2}
  \left( (ev_5^* Y + \psi_5) \cdot [\bar M_{(2,2,2,2,1)}] \right)|_{\bar M_C} =
    c_d + 3 \, [\bar M_{(2,2,2,3)}]|_{\bar M_C}.
\end {equation}

Let us first compute the virtual fundamental classes occurring in this equation.

\begin {lemma} \label {virt-fund} ~
  \begin {enumerate}
  \item The (one-dimensional) virtual fundamental class of $ \bar
    M_{(2,2,2,2,1)} $ on $ \bar M_C $ is twice the usual one.
  \item The degree of the (zero-dimensional) virtual fundamental class of $
    \bar M_{(2,2,2,3)} $ on $ \bar M_C $ is 3.
  \end {enumerate}
\end {lemma}

\begin {proof}
  We will only sketch the computations.

  (i): It is enough to do the computation at a general point of $ \bar M_C
  $. We have seen in the proof of lemma \ref {bd} that $ \bar M_{(2,2,2,2,0)} $
  is smooth of dimension 2 at a general point of $ \bar M_C $. But requiring
  multiplicity 1 at the point $ x_5 $ gives us a factor of 2 (i.e.\ $ ev_5^* Y $
  cuts out $ \bar M_C $ in $ \bar M_{(2,2,2,2,0)} $ with multiplicity 2)
  because $ x_5 $ lies on a tangency point of the stable map with $Y$.

  (ii): The only point of $ \bar M_{(2,2,2,3)} $ in $ \bar M_C $ is the stable
  map in the following picture on the left:
  \begin {center} \epsfig {file=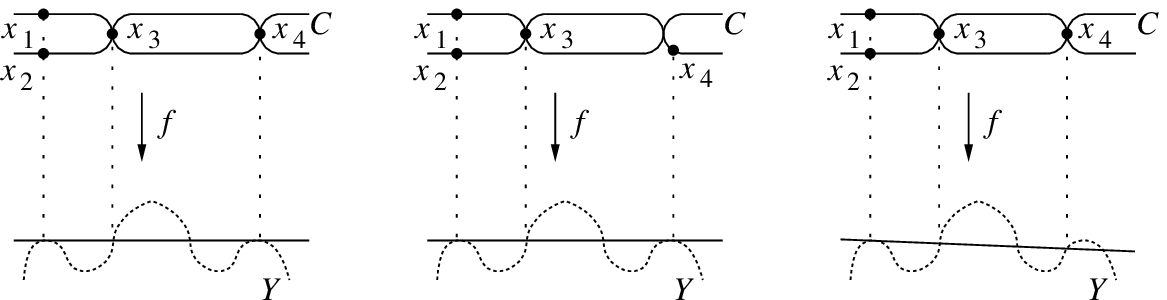} \end {center}
  Its multiplicity can be computed using remark \ref {easy-mult}. Let us study
  the space $ \bar M_{(2,2,2,2)} $ at this point first. By lemma \ref {no-def}
  every stable map in $ \bar M_{(2,2,2,2)} $ in a neighborhood of this point
  must also be a double cover of a line. It follows easily that, locally
  around this point, $ \bar M_{(2,2,2,2)} $ is reducible, with two smooth
  one-dimensional components coming together: one of them is keeping the image
  line $ f(C) $ to be the bitangent and moving the ramification point at $ x_4
  $, while keeping $ x_4 $ on an inverse image point of the bitangency point
  (picture in the middle). A local coordinate for this component is $ \epsilon
  $, where the stable map is given locally in affine coordinates as $ t \mapsto
  z=t(t+\epsilon) $, and the marked point $ x_4 $ has coordinate $ t=0 $. Now
  the equation $F$ of $Y$ vanishes on the bitangent with multiplicity 2 in $z$,
  so this equation pulled back to the curve is locally $ t^2(t+\epsilon)^2 $.
  Its $ t^2 $ coefficient vanishes to order 2 in $ \epsilon $, so this
  component contributes 2 to the virtual fundamental class of $ \bar
  M_{(2,2,2,3)} $. The other component is deforming the line in $ \PP^2 $ away
  from the bitangent, with marked points as in the picture above on the right.
  Requiring multiplicity 3 at $ x_4 $ now restricts the line back to the
  bitangent with multiplicity 1. Hence the total degree of the virtual
  fundamental class of $ \bar M_{(2,2,2,3)} $ on $ \bar M_C $ is $ 2+1=3 $.
\end {proof}

To evaluate the left hand side of equation (\ref {last2}) it is not enough to
compute the integral of $ ev_5^* Y + \psi_5 $ on $ \bar M_C $. There may also
be contributions from components of $ \bar M_{(2,2,2,2,1)} $ that just
intersect $ \bar M_C $, if the section $ \sigma $ above (whose zero locus has
class $ ev_5^* + \psi_5 $) vanishes on them at a point of $ \bar M_C $. Let us
compute these contributions.

\begin {lemma} \label {comp-4}
  The only component $Z$ of $ \bar M_{(2,2,2,2,1)} $ that meets $ \bar M_C $
  but is not $ \bar M_C $ itself corresponds to double covers of simple tangent
  lines of $Y$, with ramification and marked points as in the picture on the
  left:
  \begin {center} \epsfig {file=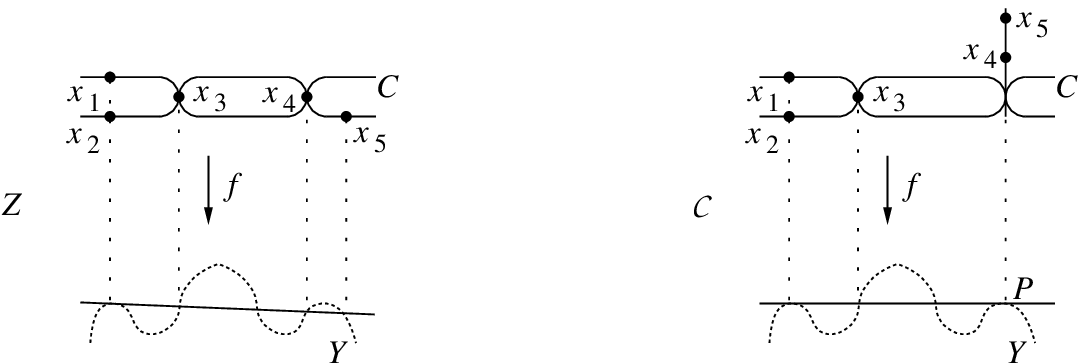} \end {center}
  Its virtual fundamental class is equal to the usual one. It intersects $ \bar
  M_C $ in the point $ \calC $ in the picture on the right. The section $
  \sigma $ on $Z$ vanishes at $ \calC $ with multiplicity 4.
\end {lemma}

\begin {proof}
  By lemma \ref {no-def} (applied to the two marked points $ x_1 $ and $ x_2 $)
  no stable map in $ \bar M_C $ can be deformed into an element in $ \bar
  M_{(2,2)} $ that is not itself a double cover of a line. So the only possible
  deformation is that we still have a double cover, however not of the
  bitangent but rather a nearby line. It is now easy to see that such a
  deformation is only possible for the stable map $ \calC \in \bar M_C $ in the
  picture above on the right, and that the deformation has to be the one in the
  picture on the left.

  The computation of the virtual fundamental class on this component $Z$ is
  completely analogous to similar calculations in previous lemmas and is
  therefore omitted. It remains to compute the order of vanishing of the
  section $ \sigma $ on $Z$ at $ \calC $. Let us forget for a moment that the
  stable map splits off a contracted rational component at this point. Let $ P
  = f(x_4) = f(x_5) $ be the bitangency point of $ \calC $. Choose local affine
  coordinates $ z_1,z_2 $ of $ \PP^2 $ around $P$ such that the local equation
  of $Y$ is $ z_2=z_1^2 + O(z_1^3) $, and $ z_2=0 $ is the bitangent. As a
  local coordinate for $Z$ around $ \calC $ we can choose $ \epsilon $, where
  the stable map is given locally around $P$ by $ t \mapsto (z_1=t^2- \frac 1 2
  \epsilon^2, z_2=\frac 1 4 \epsilon^4 + O(\epsilon^5) + t^2 \, O(\epsilon^4))
  $, and the marked points $ x_4 $ and $ x_5 $ are $ t=0 $ and $ t=\epsilon
  $. (Note that this stable map is still a double cover of a line, $ t=0 $ is a
  ramification point that maps to $Y$, and $ t=\epsilon $ another point that
  maps to $Y$.)

  Now actually the stable map splits off a rational contracted component at $
  \calC $, which corresponds to blowing up the point $ (\epsilon=0,t=0) $ in
  our family of stable maps. Hence the true local coordinates of this family of
  stable maps are the coordinates of this blow-up, i.e.\ $ t/\epsilon $ and $
  \epsilon $ instead of $t$ and $\epsilon$. So the marked points $ x_4 $ and $
  x_5 $ have coordinates $ t=0 $ and $ t/\epsilon=1 $; in particular they do
  not coincide any more at $ \epsilon=0 $.

  The vanishing of the section $ \sigma $ is the condition of tangency of $f$
  to $Y$ at $ x_5 $. So to compute its order of vanishing at $ \calC $ we
  have to look at the linear $ \xi $ coefficient of the equation of $Y$
  evaluated at the point $ t/\epsilon = 1 + \xi $, i.e.\ of
    \[ z_2 - z_1^2 + O(z_1^3)
         = \frac 1 4 \epsilon^4
           - (\frac 1 2 \epsilon^2 + \xi \epsilon^2)^2 + O(\epsilon^5)
         = - \xi \epsilon^4 - \xi^2 \epsilon^4 + O(\epsilon^5). \]
  Hence the linear $ \xi $ coefficient vanishes with multiplicity 4 at $
  \epsilon=0 $, which proves the lemma.
\end {proof}

We can now assemble the results of our local calculations to compute the number
$ c_d $.

\begin {lemma} \label {cd}
  $ c_d = -1 $ for all $d$.
\end {lemma}

\begin {proof}
  We will evaluate equation (\ref {last2})
    \[ \left( (ev_5^* Y + \psi_5) \cdot [\bar M_{(2,2,2,2,1)}]
       \right)|_{\bar M_C} =
         c_d + 3 \, [\bar M_{(2,2,2,3)}]|_{\bar M_C}. \]
  The right hand side is $ c_d + 9 $ by lemma \ref {virt-fund} (ii). The left
  hand side gets a contribution from the components of $ \bar M_{(2,2,2,2,1)} $
  that intersect $ \bar M_C $, and a contribution from $ \bar M_C \subset \bar
  M_{(2,2,2,2,1)} $ itself. The former is 4 by lemma \ref {comp-4}. The latter
  is twice the degree of $ ev_5^* Y + \psi_5 $ on $ \bar M_C $ by lemma \ref
  {virt-fund} (i). Note that the degree of $ ev_5^* Y $ on $ \bar M_C $ is
  zero, as the image point of $ x_5 $ is fixed in $ \bar M_C $.

  To compute the integral of $ \psi_5 $ on $ \bar M_C $, we will give a section
  of the cotangent line bundle $ L_5 $ and compute its zero locus. Let $z$ be a
  local coordinate around the bitangency point $ f(x_4)=f(x_5) $. Then $ f^* dz
  $ defines a section of $ L_5 $. This section vanishes only at the point where
  the moving ramification point comes to $ x_4 $ and $ x_5 $, i.e.\ at the
  point $ \calC $ in the picture of lemma \ref {comp-4} on the right. The
  computation of the order of vanishing is very similar to the calculation in
  lemma \ref {comp-4}. Ignoring the fact that $ \calC $ splits off a rational
  contracted component for $ x_4 $ and $ x_5 $, a local coordinate for $ \bar
  M_C $ around $ \calC $ is $ \epsilon $, where the stable map is given locally
  by $ t \mapsto z=t(t-\epsilon) $, the points $ x_4 $ and $ x_5 $ are $ t=0 $
  and $ t=\epsilon $, and the moving ramification point is at $ t=\frac 1 2
  \epsilon $. Now, as in the proof of the previous lemma, taking into account
  the contracted rational component of $ \calC $ means that we have to blow up
  the point $ (t=0,\epsilon=0) $, and the coordinates are actually $ t/\epsilon
  $ and $ \epsilon $. Now we see that
    \[ f^* dz = \frac {\partial}{\partial \frac t \epsilon}
                \left( \epsilon^2 \cdot \frac t \epsilon \,
                  \left( \frac t \epsilon - 1 \right) \right)
                d \frac t \epsilon
              = \epsilon^2 \left( 2 \frac t \epsilon - 1 \right)
                d \frac t \epsilon, \]
  which vanishes with multiplicity 2 in $ \epsilon $ around 0 at the point $
  x_5 $. Hence the integral of $ \psi_5 $ over $ \bar M_C $ is 2.

  Putting everything together, we get $ 4 + 2 \cdot 2 = c_d + 9 $, and therefore $
  c_d = -1 $.
\end {proof}

We can now insert the values for $ N_d $, $ b_d $ and $ c_d $ from proposition
\ref {Nd}, lemma \ref {bd}, and lemma \ref {cd}, respectively, into the
equation from lemma \ref {nd-cor}, and get the following final result.

\begin {corollary} \label {nd}
  For $ d \ge 5 $, the enumerative number of conics 5-fold tangent to $Y$ is
    \[ n_d = \frac 1 {5!} \, d \, (d-3) \, (d-4) \,
             (d^7+12d^6-18d^5-540d^4+251d^3+5712d^2-1458d-14580). \]
\end {corollary}

\begin {remark}
  The first few values of $ n_d $ are given in the following table.
    \[ \begin {array}{|c|c|c|c|c|c|c|} \hline
         d & 5 & 6 & 7 & 8 & 9 & 10 \\ \hline
         n_d & 2015 & 70956 & 684222 & 3878736 & 16137873 &  54575640 \\ \hline
       \end {array} \]
\end {remark}


\section {Application to rational curves on K3 surfaces} \label {sec-k3}

Let $X$ be a K3 surface, and let $ \beta \subset H_2 (X,\ZZ) $ be the class of a
holomorphic curve in $X$. The moduli space of stable maps to $X$ of class $
\beta $ has virtual dimension $ -1 $, so there is no corresponding
Gromov-Witten invariant of $X$. However, we have chosen $X$ such that it
contains rational curves of class $ \beta $, and we would like to count
them. The reason for the mismatch in the virtual dimension is that, in the
space of all K3 surfaces, only a 1-codimensional subset of K3 surfaces contains
curves in the class $ \beta $ at all. So, if $ \calX $ is a general
1-dimensional family of K3 surfaces with $X$ as central fiber, and $ \tilde
\beta \in H_2 (\calX,\ZZ) $ is the class induced by $ \beta $ via the inclusion
$ X \subset \calX $, then the only rational curves in $ \calX $ of class $
\tilde \beta $ are in fact curves of class $ \beta $ in $X$. But now the virtual
dimension of the space of stable maps to $ \calX $ of class $ \beta $ is 0,
hence there is a corresponding Gromov-Witten invariant. This invariant counts
curves in $ \calX $ of class $ \tilde \beta $ and therefore curves in $X$ of
class $ \beta $; so we would like to call this number $ n_\beta$ ``the number
of rational curves in $X$ of class $ \beta $''. A rigorous definition of the
invariant $ n_\beta $ of $X$ along these lines has been given in \cite
{BL}. The number $ n_\beta $ does not depend on a family $ \calX $ chosen to
define it.

The numbers $ n_\beta $ have been computed in various papers (\cite {YZ}, \cite
{B}, \cite {Go}, \cite {BL}) under the assumption that the class $ \beta $ is
primitive, i.e.\ not a non-trivial multiple of a smaller integral homology
class. The result is that $ n_\beta $ is equal to the $ q^d $ coefficient in
the series
\begin {align*}
  G(q) = \prod_{i>0} \frac 1 {(1-q^i)^{24}}
      &= 1 + 24q + 324q^2 + 3200q^3
           + 25650q^4 + 176256q^5 + \cdots \\
      &=: \sum_{d \ge 0} G_d q^d.
\end {align*}
where $ d = \frac 1 2 \beta^2 + 1 $. It is not known yet what the numbers are
if $ \beta $ is not primitive.

The results of this paper allow us to compute the number $ n_\beta $ explicitly
in a case where $ \beta $ is not primitive. Let $ Y \subset \PP^2 $ be a
general sextic curve, and let $ \pi: X \to \PP^2 $ be the double cover of $
\PP^2 $ branched along $Y$. It is well-known that $X$ is a K3 surface. Let us
start by considering curves on $X$ in the (primitive) class $ \beta = \pi^*
\ell $, where $ \ell $ is the class of a line. The pull-back of a general line
in $ \PP^2 $ will be a 2:1 cover of $ \PP^1 $, branched along the 6
intersection points of $Y$ with the line, hence it is a curve of genus 2. We
get a rational curve (with 2 nodes) on $Y$ as a pull-back of a line if and only
if the line is a bitangent of $Y$: the pull-back is then again a 2:1 cover, but
with 2 nodes (the bitangency points), and only 2 ramification points (the
remaining 2 intersection points of $Y$ with the line). So we see that $
n_{\pi^* \ell} $ has to be the number of bitangents of $Y$. In fact, this
number is 324 (see e.g.\ \cite {H} exercise IV.2.3 f), which is equal to $ G_2
$ (note that $ \frac 1 2 (\pi^* \ell)^2 + 1 = 2 $).

Now let us consider rational curves in $X$ of class $ 2 \pi^* \ell $, i.e.\
pull-backs of conics --- this class is not primitive any more. The pull-back of
a general conic will be a 2:1 cover of the conic ramified at 12 points, so it
is a curve of genus 5. We can get rational curves in the following ways:
\begin {enumerate}
\item Pull-backs of (smooth) conics that are 5-fold tangent to $Y$. These will
  be 2:1 covers of the conic, with 5 nodes over the tangency points, and only 2
  ramification points (the remaining 2 intersection points of $Y$ with the
  conic). By corollary \ref {nd}, there are $ 70956 $ such curves.
\item Pull-backs of unions of two distinct lines: these give a rational curve
  only if the pull-backs of both lines are rational, i.e.\ they are both
  bitangents. The pull-back is then a union of two (2-nodal) rational curves on
  $X$ that intersect in 2 points. To make this into a rational stable map we
  can glue these two components at either intersection point. Hence there are
  $ 2 \cdot \binom {324}{2} = 104652 $ such stable maps.
\item Double covers of pull-backs of a line, necessarily again of a bitangent.
  The pull-back of such a bitangent is a 2-nodal rational curve $C$. There are
  two possible ways of double covers of such a curve:
  \begin {enumerate}
  \item Double covers that factor through the normalization of $C$. The space
    of these curves is the same as that of double covers of a smooth rational
    curve; it has dimension 2.
  \item Double covers that do not factor through the normalization. They have
    two components that are both mapped to $C$ with degree 1, and glued over
    one of the nodes of $C$ in such a way that, locally around this node, the
    morphism of the stable map is an isomorphism onto $C$. There are $ 2 \cdot
    324 = 648 $ such curves.
  \end {enumerate}
\end {enumerate}
Adding up just the numbers from (i), (ii), and (iii b), we get
  \[ 70956 + 104652 + 648 = 176256, \]
which is exactly $ G_5 $ (and $ \frac 1 2 (2 \pi^* \ell ) + 1 = 5 $). So we see
that, for our non-primitive class $ \beta $, the corresponding invariant from
the series $ G(q) $ does give the correct number, \emph {except} for a
correction (iii a) for double covers of curves of class $ \frac 1 2 \beta $
that factor through the normalization of these curves. Let us compute what this
correction term is.

\begin {lemma}
  With notations as above, the double covers of type (iii a) of the pull-back
  of a bitangent contribute $ \frac 1 8 $ to the invariant $ n_{2\pi^* \ell} $.
\end {lemma}

\begin {proof}
  Let $ D \cong \PP^1 $ be the normalization of the nodal rational curve in
  $X$. The moduli space of the double covers that factor through the
  normalization is then just $ \bar M_{0,0} (D,2) $, which has dimension 2. As
  the normal bundle of the (local) immersion $ D \to X $ is $ \calO(-2) $, the
  rank-3 obstruction bundle for the corresponding Gromov-Witten invariant would
  be $ R^1\pi_* f^* \calO(-2) $, where $ \pi: \bar M_{0,1} (\PP^1,2) \to \bar
  M_{0,0} (\PP^1,2) $ is the forgetful map and $ f: \bar M_{0,1} (\PP^1,2) \to
  \PP^1 $ the evaluation.

  As explained above, the K3-invariants of $X$ are defined as the ordinary
  Gromov-Witten invariants of a 1-dimensional family $ \calX $ of K3 surfaces
  in which $X$ is the only surface that contains rational curves in the given
  homology class. This means that the obstruction bundle for the K3 invariants
  is obtained from the usual Gromov-Witten obstruction bundle by taking the
  quotient with $ \pi_* f^* N_{X/\calX} = \pi_* f^* \calO = \calO $. So the
  integral that we want to compute is
    \[ c_{top} (R^1 \pi_* f^* \calO(-2) / \calO) \cdot
         [\bar M_{0,0} (\PP^1,2)]. \]
  This is easily done: from the two exact sequences on $ \bar M_{0,1} (\PP^2,2)
  $
  \begin {gather*}
    0 \to f^* \calO(-1) \to f^* \calO \to f^* \calO_P \to 0 \\
    0 \to f^* \calO(-2) \to f^* \calO(-1) \to f^* \calO_P \to 0
  \end {gather*}
  (where $ P \in \PP^1 $ is a point) we get the exact sequences of vector
  bundles on $ \bar M_{0,0} (\PP^2,2) $
  \begin {gather*}
    0 \to \calO \to \pi_* f^* \calO_P \to R^1 \pi_* f^* \calO(-1) \to 0 \\
    0 \to \pi_* f^* \calO_P \to R^1 \pi_* f^* \calO(-2)
      \to R^1 \pi_* f^* \calO(-1) \to 0
  \end {gather*}
  from which it follows that
    \[ c_2 (R^1 \pi_* f^* \calO(-2) / \calO) =
       c_2 (R^1 \pi_* f^* (\calO(-1) \oplus \calO(-1))), \]
  i.e.\ the contribution of the double covers under consideration is the same
  as the double cover contribution for rational curves on Calabi-Yau
  threefolds with balanced normal bundle. This contribution is well-known to be
  $ \frac 1 8 $, see e.g.\ \cite {G2} example 6.3.
\end {proof}

So we see that
  \[ n_{2\pi^* \ell} = G_5 + \frac 1 8 \cdot G_2. \]
We conjecture that this pattern continues, i.e.\ that the numbers $ n_\beta $
receive multiple cover corrections similarly to the case of Gromov-Witten
invariants of Calabi-Yau threefolds:
  \[ n_\beta = \sum_k \frac 1 {k^3} \cdot G_{\frac 1 2 (\frac \beta k)^2 + 1}, \]
where the sum is taken over all $k>0$ such that $ \frac \beta k $ is an
integral homology class.


\end {document}